\documentclass[11pt]{amsart}
\usepackage{amsmath}
\usepackage{amsfonts}

\usepackage{amscd}
\usepackage{amsthm}
\usepackage{amssymb} \usepackage{latexsym}
\usepackage{eufrak}
\usepackage{euscript}
\usepackage{epsfig}
\usepackage{graphics}
\usepackage{array}
\usepackage{enumerate}
\usepackage{dsfont}
\usepackage{color}
\usepackage{wasysym}
\usepackage{hyperref}
\usepackage{pdfsync}

\newcommand{\bel}[1]{\begin{equation}\label{#1}}

\newcommand{\be}{\begin{equation}}

\newcommand{\ba}{\begin{eqnarray}}
\newcommand{\ea}{\end{eqnarray}}

\newcommand{\qe}{\end{equation}}
\newcommand{\R}{{\mathbb R}}
\newcommand{\N}{{\mathbb N}}
\newcommand{\Z}{{\mathbb Z}}

\newcommand{\CCD}{\mathrm{CD}}

\newcommand{\DD}{\mathrm{D}}
\newcommand{\gir}{\mathrm{Gir}}
\newcommand{\str}{\mathrm{Star}}
\newcommand{\dd}{\mathrm{d}}

\newcommand{\Hmm}[1]{\leavevmode{\marginpar{\tiny%
$\hbox to 0mm{\hspace*{-0.5mm}$\leftarrow$\hss}%
\vcenter{\vrule depth 0.1mm height 0.1mm width \the\marginparwidth}%
\hbox to
0mm{\hss$\rightarrow$\hspace*{-0.5mm}}$\\\relax\raggedright #1}}}

\theoremstyle{theorem}
\newtheorem{thm}{Theorem}[section]
\newtheorem{prop}[thm]{Proposition}
\theoremstyle{example}

\theoremstyle{corollary}
\newtheorem{coro}[thm]{Corollary}
\theoremstyle{lemma}
\newtheorem{lemma}[thm]{Lemma}
\theoremstyle{definition}
\newtheorem{defi}[thm]{Definition}
\theoremstyle{proof}
\theoremstyle{remark}

\begin{document}

\title[Large girth and curvature dimension condition]{Graphs with large girth and nonnegative curvature dimension condition}
\author{Bobo Hua}
\email{bobohua@fudan.edu.cn}
\address{School of Mathematical Sciences, LMNS, Fudan University, Shanghai 200433, China
}

\author{Yong Lin}
\email{linyong01@ruc.edu.cn}
\address{Department of Mathematics,Information School,
Renmin University of China,
Beijing 100872, China}

\begin{abstract} In this paper, we classify unweighted graphs satisfying the curvature dimension condition $\CCD (0,\infty)$ whose girth are at least five.
\end{abstract}
\maketitle%\tableofcontents

\section{introduction}

In Riemannian geometry, there are various geometric curvature notions, such as sectional curvature, Ricci curvature and scalar curvature, derived from the Riemann curvature tensor. Of particular interest,  curvature bounds usually impose many topological and geometric constraints for underlying manifolds. Even in the non-smooth setting, there are generalizations of curvature bounds, e.g. sectional curvature on Alexandrov spaces, see \cite{BuragoGromovPerelman92,Burago01}, and Ricci curvature on metric measure spaces \cite{LottVillani09,Sturm06a,Sturm06b}, from which many geometric consequences can be derived accordingly. 

Many authors attempted to define appropriate curvature conditions on discrete metric spaces, e.g. graphs, in order to resemble some geometric properties of Riemannian curvature bounds. 

One is so-called combinatorial curvature introduced by \cite{Stone:1976wz,Gromov87,Ishida90}.  The idea is to properly embed a graph into a Riemannian manifold, in particular a surface, and to define the curvature bound of the graph from that of the ambient space.  In this way, one can derive some global geometric properties of the graph via the embedding, see \cite{Zuk:1997jc,WOESS:1998jy,Higuchi:2001ft,Baues:2001bs,BauesPeyerimhoff06,DeVos:2007fq,Chen:2008ba,Chen:2008dm,Keller:2008vo,Keller:2011gg,KPF14,Hua:2015gl}.

Ollivier \cite{Ollivier09} used $L^1$-Wasserstein distance for the space of probability measures on graphs to define a curvature notion mimicking the Ricci curvature on manifolds. Interesting results can be obtained from the optimal transport strategy, see \cite{BauerJostLiu12,OllivierVillani12,Paeng12,JostLiu14,BhaMuk15}. Lin, Lu and Yau \cite{LinLuYau11} modified Ollivier's definition and \cite{LinLuYau14} gave a classification of Ricci flat graphs with girth at least five.  Maas \cite{Maas11} identified the heat flow and the gradient flow of the Boltzmann-Shannon type entropy by introducing a Riemannian structure on the space of probability measures on graphs. Erbar and Maas \cite{ErbarMaas12} defined the generalized Ricci curvature via the convexity of the entropy functional and derived many functional inequalities under this curvature assumption.

From a different strategy, one can define curvature dimension conditions via the so-called $\Gamma$-calculus for general Markov semigroups, where $\Gamma$ is the ``carr\'e du champ" operator, see \cite[Definition~1.4.2]{BakryGentilLedoux}. In particular, the curvature bound is defined via a Bochner type inequality using the iterated $\Gamma$ operator, denoted by $\Gamma_2,$ see Definition~\ref{d:carre du} in this paper. For the diffusion semigroup, curvature dimension conditions were initiated in Bakry and \'Emery \cite{BakryEmery85}, and for the non-diffusion case, e.g. graphs, 
introduced by Lin and Yau \cite{LinYau10}. Later, variants of curvature dimension conditions were introduced to obtain important analytic results, see e.g. \cite{HLLY14,HLLY14,Muench14,Bauer:2015vn,LinLiu15,Hua:2015tv,GongLin15,Muench15,FathiShu15,KKRT16,CLP16}.

We introduce the setting of graphs and refer to Section~\ref{s:graphs} for details.
Let $(V,E)$ be an undirected, connected, locally finite simple graph with the set of vertices $V$ and the set of edges $E.$ Without loss of generality, we exclude the trivial graph consisting of a single vertex. Two vertices $x,y$ are called neighbors if $\{x,y\}\in E,$ denoted by $x\sim y.$ The combinatorial degree of a vertex $x\in V$ is the number of its neighbors, denoted by $\dd_x.$ We assign a weight $m_x$ to each vertex $x$ and a weight $\mu_{xy}$ to each edge $\{x,y\},$ and refer to the quadruple $G=(V,E,m,\mu)$ as a \emph{weighted graph}. The graph $G$ is called \emph{unweighted} if $\mu\equiv 1$ on $E.$ For any $x\in V,$ we denote $\mu_x:=\sum_{y\sim x}\mu_{xy}.$ 

We are mostly interested in functions defined on $V,$ and denote by $C(V)$ the set of all such functions.
For any weighted graph $G$, there is an associated \emph{Laplacian} operator, $\Delta:C(V)\to C(V),$ defined as
\begin{equation}\label{defLaplacian}\Delta f(x)=\frac{1}{m_x}\sum_{y\sim x}\mu_{xy}(f(y)-f(x)),\quad f\in C(V), x\in V.\end{equation} One can see that the weights $\mu$ and $m$ play the essential role in the definition of Laplacian.
Given the weight $\mu$ on $E,$ typical choices of $m$ are of interest:
\begin{itemize}\item In case of $m_x=\mu_x$ for any $x\in V,$ we call the associated Laplacian the \emph{normalized} Laplacian.
\item In case of $m\equiv1$ on $V,$ the Laplacian is called \emph{physical} (or combinatorial) Laplacian.\end{itemize} Moreover, if the graph is unweighted, the corresponding Laplacian is called unweighted normalized (i.e. $\mu\equiv 1$ on $E$ and $m\equiv \mu$ on $V$) or unweighted physical Laplacian (i.e. $\mu\equiv 1$ on $E$ and $m\equiv 1$ on $V$) respectively. For simplicity, we also call the graph \emph{unweighted normalized} 
or \emph{unweighted physical} graph accordingly.

We denote by $\ell^p(V,m)$ or simply $\ell^p_m,$ the space of $\ell^p$ summable functions on the discrete measure space $(V,m)$ and by $\|\cdot\|_{\ell^p_m}$ the $\ell^p$ norm of a function.
 Define the weighted vertex degree $\mathrm{D}:V\to[0,\infty)$ by
\begin{align*}
    \mathrm{D}_x=\frac{1}{m_x}\sum_{y\sim x}\mu_{xy},\qquad x\in V.
\end{align*} It is well known, see e.g. \cite{KellerLenz12}, that the Laplacian associated with the graph $G$ is a bounded operator on $\ell^2_m$
if and only if $\sup_{x\in V} \mathrm{D}_x<\infty.$

The curvature dimension condition $\CCD(K,n),$ for $K\in \R$ and $n\in (0,\infty],$ on graphs was introduced by \cite{LinYau10}, which serves as the combination of a lower bound $K$ for the Ricci curvature and an upper bound $n$ for the dimension, see Definition \ref{d:curvature dimension}.   To verify the $\CCD(K,n)$ condition, we adopt the following crucial identity for general Laplacians, analogous to the Bochner identity on Riemannian manifolds, which was first proved in   
\cite{LinYau10}, see also \cite{Mali13}, for normalized Laplacians. 
\begin{prop}\label{propboch}For any function $f$ and $x\in V,$
\begin{equation}\label{bochner}\Gamma_2(f)(x)=\frac14 |D^2 f|^2(x)+\frac12(\Delta f (x))^2-\frac14\sum_{y\sim x}\frac{\mu_{xy}}{m_x}(\DD_x+\DD_y)|f(y)-f(x)|^2,\end{equation} where 
$$|D^2 f|^2(x):=\sum_{\substack{y,z\in V\\ y\sim x, z\sim y}}\frac{\mu_{xy}\mu_{yz}}{m_xm_y}|f(x)-2f(y)+f(z)|^2.$$
\end{prop} Note that $|D^2 f|^2$ above is a discrete analogue of the squared norm of the Hessian of a function $f$ in the Riemannian setting.

\begin{figure}[htbp]
\begin{center}
\includegraphics[width=\linewidth]{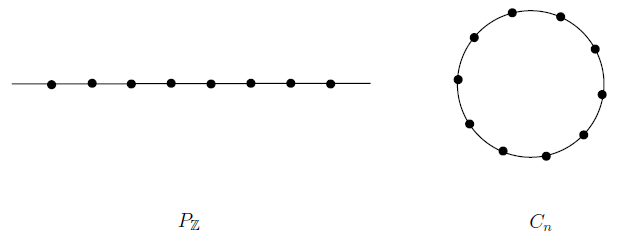}
\caption{Theorem~\ref{main-normalized1}}
\label{fig1}
\end{center}
\end{figure}

\begin{figure}[htbp]
\begin{center}
\includegraphics[width=\linewidth]{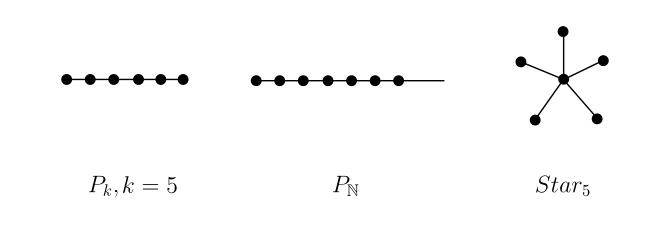}
\caption{Theorem~\ref{main-normalized2}}
\label{fig2}
\end{center}
\end{figure}

\begin{figure}[htbp]
\begin{center}
\includegraphics[width=\linewidth]{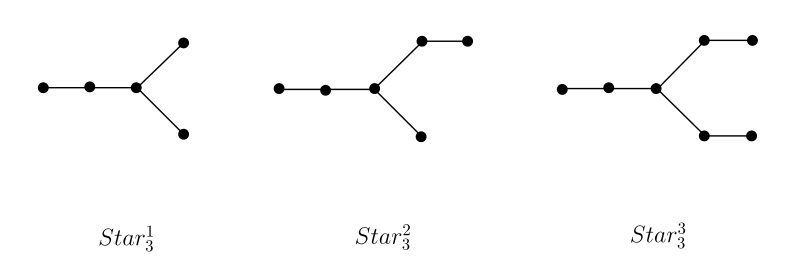}
\caption{Theorem~\ref{main-normalized2}}
\label{fig3}
\end{center}
\end{figure}

The girth of a vertex $x,$ denoted by $\gir(x),$ is defined as the minimal length of cycles passing through $x,$ and the girth of a graph is the minimal girth of vertices, see Definition~\ref{girth1}. Inspired by the work \cite{LinLuYau14}, we classify the unweighted graphs with large girth and satisfying the $\CCD(0,\infty)$ condition. By definition, the curvature condition at a vertex is determined by the local structure, in particular, the ball of radius two centered at the vertex, denoted by $B_2.$ The key observation is that if the girth of a vertex is large, $B_2$ is essentially a tree, see Proposition~\ref{tree}, which is intuitively non-positively curved. By using the Bochner type identity \eqref{bochner}, one obtains the sufficient and necessary condition for $\CCD(0,\infty)$ in that case, see Corollary~\ref{equCD0}, which yields the following classifications.
%We say $x,y\in V$ are neighbors, denoted by $x\sim y,$ if $(x,y)\in E.$ The graph is called \emph{locally finite} if each vertex has finitely many neighbors. In this paper, we only consider locally finite graphs.  We denote by $$C_0(V):=\{f:V\to\R|\ \{x\in V| f(x)\neq 0\}\ \mathrm{is\ of\ finite\ cardinality}\}$$ the set of finitely supported functions on $V$ and by $\ell^p(V,m),$ $p\in [1,\infty],$ the $\ell^p$ spaces of functions on $V$ with respect to the measure $m.$
\begin{thm}\label{main-normalized1} Let $(V,E,m,\mu)$ be an unweighted normalized graph satisfying the $\CCD(0,\infty)$ condition and $\inf_{x\in V}\dd_x\geq 2.$ If for some $x_0\in V,$ $\gir(x_0)\geq 5,$ then the graph is either the infinite line $P_{\Z}$ or the cycle graphs $C_n$ for $n\geq 5,$ see Figure~\ref{fig1}.
\end{thm} 
It is remarkable that if all vertex degrees are at least two, to derive the classification we only assume the girth of a vertex in the graph is large. For the general case below, a stronger assumption that the girth of the whole graph is large is needed.
\begin{thm}\label{main-normalized2} Let $(V,E,m,\mu)$ be an unweighted normalized graph satisfying the $\CCD(0,\infty)$ condition whose girth is at least $5.$ Then the graph is one of the following:
\begin{enumerate}[(a)]
\item The path graphs $P_k$ $(k\geq 1),$ the cycle graphs $C_n$ $(n\geq 5),$ the infinite line $P_{\Z},$ or the infinite half line $P_{\N},$ see Figure~\ref{fig2}.
\item The star graphs $\str_n$ $(n\geq 3),$ or $\str_{3}^{i}$ $(1\leq i\leq 3),$
\end{enumerate} where $\str_3^i$ is the $3$-star graph with $i$ edges added, $1\leq i\leq 3,$ see Figure~\ref{fig3}. 
\end{thm} 

For physical Laplacians, we also obtain the classification results, see Section~\ref{s:physical}. Note that, similar results for physical Laplacians have been obtained in Cushing, Liu and Peyerimhoff \cite[Corollary~6.9]{CLP16}.

The organization of the paper is as follows:
In next section, we introduce the definitions for graphs, $\Gamma$-calculus, and criteria for curvature dimension conditions for graphs with large girth. In Section~\ref{s:normal}, we study normalized graphs and prove the classification results, Theorem~\ref{main-normalized1} and \ref{main-normalized2}. The last section is devoted to physical Laplacians.

\section{Graphs}\label{s:graphs}
\subsection{Combinatorial and weighted graphs}
Let $(V,E)$ be a (finite or infinite) undirected graph with the set of vertices $V$ and the set of edges $E,$ i.e. two-elements subsets in $V.$ The graph is called simple if there is no self-loops and multiple edges. The graph is called locally finite, if the combinatorial degree $\dd_x<\infty$ for any $x\in V.$ We say a vertex $x$ is a pending vertex if $\dd_x=1.$ For any subsets $A,B\subset V,$ we denote by $E(A,B):=\{\{x,y\}\in E: x\in A, y\in B\}$ the set of edges between $A$ and $B.$ For vertices $x$ and $y,$ a \emph{walk} from $x$ to $y$ is a sequence of vertices $\{x_i\}_{i=0}^k$ such that $$x=x_0\sim x_1\sim \cdots\sim x_k=y,$$ where $k$ is called the length of the walk.  A graph is said to be connected if for any $x,y\in V$ there is a walk from $x$ to $y.$ The minimal length of walks from $x$ to $y$ is called the (combinatorial) distance between them, denoted by $d(x,y).$ In this paper, we only consider undirected, connected, locally finite simple graphs. 

A \emph{cycle} of length $k>2$ is a walk $\{x_i\}_{i=0}^k$ such that $x_0=x_k$ and $x_i\neq x_j$ for all $i\neq j, 0\leq i,j\leq k-1.$ A graph is called a tree if it contains no cycles.

\begin{defi}\label{girth1} The girth of a vertex $x$ in $(V,E),$ denoted by $\gir(x)$, is defined to the minimal length of cycles passing through $x.$ (If there is no cycle passing through $x,$ define $\gir(x)=\infty.$) The girth of a graph is defined as $\inf_{x\in V}\gir(x).$
\end{defi}

For any $x\in V,$ $r\in \N_0,$ we denote by $B_r(x):=\{y\in V: d(y,x)\leq r\}$ the ball of radius $r$ centered at $x,$ and by $S_r(x):=\{y\in V: d(y,x)= r\}$ the corresponding sphere. For our purposes, we define a graph, denoted by $\widehat{B_2}(x),$ consisting of the set of vertices in $B_2(x)$ and the set of edges $\{\{x,y\}\in E: x\in B_1(x)\ \mathrm{or}\ y\in B_1(y)\}.$ That is, $\widehat{B_2}(x)$ is obtained by removing edges in $E(S_2(x),S_2(x))$ from the induced subgraph $B_2(x).$ 
The following proposition is elementary and useful.
\begin{prop}\label{tree} For a graph $(V,E)$ and $x\in V,$ $\widehat{B_2}(x)$ is a tree if and only if $\gir(x)\geq 5.$
\end{prop}
\begin{proof} $\Longrightarrow:$ Suppose $\widehat{B_2}(x)$ is a tree and $\gir(x)\leq 4.$ Let $C=\{x_i\}$ be a cycle of minimal length passing through $x$ of length $\leq 4.$ Then $C$ is a cycle in $\widehat{B_2}(x)$ which contradicts to that $\widehat{B_2}(x)$ is a tree.

$\Longleftarrow:$ Conversely, suppose that $\gir(x)\geq 5$ and $\widehat{B_2}(x)$ is not a tree, then there is a cycle $C=\{x_i\}$ in $\widehat{B_2}(x).$ We divide it into two cases:

\emph{Case 1.} If the cycle $C$ contains no vertices in $S_2(x),$ then it is included in $B_1(x).$ The cycle $C$ passes through an edge $e=(x_j,x_{j+1})$ in $E(S_1(x),S_1(x)).$ (Otherwise it will be contained in the graph $\widehat{B_1}(x)$ obtained by removing edges $E(S_1(x),S_1(x))$ from the induced subgraph $B_1(x)$ which is a tree.  A contradiction.) Now we have a cycle $\{x,x_j,x_{j+1}\}$ of length $3$ which contradicts to $\gir(x)\geq 5.$

\emph{Case 2.} $S_2(x)\cap C\neq \emptyset.$ Let $x_k\in S_2(x)$ for some $k$ and, without loss of generality, denote $x_{k-1}\sim x_k\sim x_{k+1}$ in $C.$ Since in $\widehat{B_2}(x)$ there is no edges connecting vertices in $S_2(x),$ the consecutive neighbors of $x_k$ in $C$, $x_{k-1}$ and $x_{k+1},$ are hence contained in $S_1(x).$ Thus, $\{x,x_{k-1},x_k,x_{k+1}\}$ is a cycle of length $4.$ A contradiction.

\end{proof} 

As mentioned in the introduction, for a combinatorial graph $(V,E),$ we assign weights on the set of vertices $V$ and edges $E$ respectively, $m:V\to(0,\infty)$ and $\mu:E\to(0,\infty),$ to obtain a weighted graph $G=(V,E,m,\mu).$ In the following we always write $G$ abbreviately for a weighted graph. For convenience, we extend the function $\mu$ on $E$ to the total set $V\times V,$ $\mu:V\times V\to [0,\infty),$ by 
%$(x,y)\mapsto \mu_{xy}$ if $x\sim y,$ and $\mu_{xy}=0$ for any $x\not\sim y.$ 
\[
(x,y)\mapsto \left\{
\begin{array}{ccc}
  \mu_{xy},& \mathrm{for}\ x\sim y,  \\
 0,&\mathrm{for}\ x\not\sim y.\\
\end{array}
\right.
\] So that we may write for $x\in V,$
$$\sum_{y}\mu_{xy}f(x,y)=\sum_{y\sim x}\mu_{xy}f(x,y)$$ in the following context.
The Laplacian of a weighted graph $G$ is defined as in \eqref{defLaplacian} which can be identified with the generator of a standard Dirichlet form associated to the weighted graph $G,$ see \cite{KellerLenz12}. 

\subsection{Gamma calculus}\label{subsection2.2}
We introduce the $\Gamma$-calculus and curvature dimension conditions on graphs following \cite{LinYau10}.%,Bauer13}.

Given $f: V\to\R$ and $x,y\in V,$ we denote by $\nabla_{xy}f:= f(y)-f(x)$ the difference of the function $f$ on the vertices $x$ and $y.$ First we define two natural bilinear forms associated to the Laplacian $\Delta$. 
\begin{defi}\label{d:carre du}
  The gradient form $\Gamma,$ called the ``carr\'e du champ" operator, is defined by, for $f,g\in C(V)$ and $x\in V$,
  \begin{eqnarray*}\Gamma(f,g)(x)&=&\frac12(\Delta(fg)-f\Delta g-g\Delta f)(x)\\&=&\frac{1}{2m(x)}
  \sum_{y}\mu_{xy}\nabla_{xy}f\nabla_{xy}g.\end{eqnarray*} For simplicity, we write $\Gamma(f):=\Gamma(f,f).$ Moreover, the iterated gradient form, denoted by $\Gamma_2$, is defined as $$\Gamma_2(f,g)=\frac12(\Delta\Gamma(f,g)-\Gamma(f,\Delta g)-\Gamma(g,\Delta f)),$$ and we write $\Gamma_2(f):=\Gamma_2(f,f)=\frac{1}{2}\Delta \Gamma(f)-\Gamma(f,\Delta f).$
\end{defi}

%The Cauchy-Schwarz inequality implies that \begin{equation}\label{eq:gamma}\Gamma(f,g)\leq \sqrt{\Gamma(f)\Gamma(g)}\leq \frac12 (\Gamma(f)+\Gamma(g)).\end{equation} In addition, one can easily see that $Q^{(N)}(f)=\|\Gamma(f)\|_{\ell^1_m}.$

Now we prove the Bochner type identity on graphs.
\begin{proof}[Proof of Proposition~\ref{propboch}]
For any function $f$ and $x\in V,$
\begin{eqnarray*}\Delta\Gamma(f)(x)&=&\sum_{y}\frac{\mu_{xy}}{m_x}\Gamma(f)(y)-\DD_x\Gamma(f)(x)
=\sum_{y,z}\frac{\mu_{xy}\mu_{yz}}{2m_xm_y}(\nabla_{yz}f)^2-\DD_x\Gamma(f)(x)\\
&=&\sum_{y,z}\frac{\mu_{xy}\mu_{yz}}{2m_xm_y}\left[(\nabla_{yz}f)^2-(\nabla_{xy}f)^2\right]+\sum_y\frac{\mu_{xy}}{2m_x}(\DD_y-\DD_x)(\nabla_{xy}f)^2\\
&=& (I)+\sum_y\frac{\mu_{xy}}{2m_x}(\DD_y-\DD_x)(\nabla_{xy}f)^2.
\end{eqnarray*} For the first term on the right hand side of the equation, we write
\begin{eqnarray*}(I)&=&\sum_{y,z}\frac{\mu_{xy}\mu_{yz}}{2m_xm_y}
\left\{(\nabla_{yz}f-\nabla_{xy}f)^2+2\nabla_{xy}f(\nabla_{yz}f-\nabla_{xy}f)\right\}\\&=&\frac12|D^2f|^2(x)+\sum_y\frac{\mu_{xy}}{m_x}\nabla_{xy}f\Delta f(y)-\sum_{y}\frac{\mu_{xy}}{m_x}\DD_y(\nabla_{xy}f)^2\\
&=&\frac12|D^2f|^2(x)+2\Gamma(f,\Delta f)+(\Delta f(x))^2-\sum_{y}\frac{\mu_{xy}}{m_x}\DD_y(\nabla_{xy}f)^2.
\end{eqnarray*} Combining the above equations, we prove the proposition.
\end{proof}

Now we can introduce curvature dimension conditions on graphs.
\begin{defi}\label{d:curvature dimension}
 Let $K\in \R, n\in (0,\infty].$ We say a graph $G$ satisfies the $\CCD(K,n)$ condition at $x\in V,$ denoted by $\CCD(K,n,x),$ if for any $f\in C(V)$,
  \begin{equation}\label{CDcondition}\Gamma_2(f)(x)\geq \frac1n(\Delta f(x))^2+K\Gamma (f)(x).\end{equation} A graph is said to satisfy $\CCD(K,n)$ condition if the above inequality holds for all $x\in V.$ 
\end{defi}

\subsection{Criteria for curvature dimension conditions} By Proposition~\ref{propboch}, we have the following criterion for the curvature dimension condition of a vertex with large girth.
\begin{thm}\label{main1} If $\gir(x)\geq 5,$ then $\CCD(K,n, x)$ holds if and only if for any $f\in C(V),$ 
\begin{equation}\label{equivalent}(1-\frac{2}{n})(\Delta f(x))^2\geq \sum_y \frac{\mu_{xy}}{m_x}\left(\frac{\DD_x+\DD_y}{2}-\frac{2\mu_{xy}}{m_y}+K\right)(f(y)-f(x))^2.\end{equation}
\end{thm}
\begin{proof} Since the terms $\Delta f, \Gamma(f)$ and $\Gamma_2(f)$ are all invariant by adding a constant to $f$,
it suffices to check the curvature conditions at $x\in V$ for functions $f$ satisfying $f(x)=0.$ Note that the right hand side of \eqref{CDcondition} only depends on the values of $f$ on $B_1(x).$ Set $W_f:=\{g\in C(V): g|_{B_1(x)}=f|_{B_1(x)}\}.$ It suffices to prove the following 
\begin{equation}\label{eq1}\inf_{g\in W_f} \Gamma_2(g)=\sum_{y}\frac{\mu_{xy}^2}{m_xm_y}f(y)^2+\frac12 (\Delta f(x))^2-\frac14\sum_{y}\frac{\mu_{xy}}{m_x}(\DD_x+\DD_y)f(y)^2.\end{equation}
By the formula in \eqref{bochner}, it suffices to minimize
$|D^2 g|^2(x)$ under the same constraints. Note that $\widehat{B_2}(x)$ is a tree by Proposition~\ref{tree}, i.e. for any $z\in S_2(x)$ there is a unique path from $z$ to $x.$ For $g\in W_f,$
\begin{eqnarray*}|D^2 g|^2(x)&=&\sum_{y\in V}\frac{\mu_{xy}}{m_xm_y}\sum_{z\in S_2(x)\cup \{x\}}\mu_{yz}|g(z)-2f(y)|^2\\
&=&\sum_{y\in V}\frac{\mu_{xy}}{m_xm_y}\left(\sum_{z\in S_2(x)}\mu_{yz}|g(z)-2f(y)|^2+\mu_{yx}|2f(y)|^2\right)
.\end{eqnarray*} The first equality follows from the fact that the nontrivial terms in the summation are all in the form $x\sim y\sim z,$ and hence $z\in S_2(x)\cup \{x\}.$ Then it is easy to see that the infimum over $g\in W_f$ is attained by setting $g(z)=2f(y)$ for any $z\in S_2(x)$ where $y$ is the unique vertex in $S_1(x)$ such that $x\sim y\sim z.$ This proves \eqref{eq1} and hence the theorem.
\end{proof}

For the curvature conditions at $x\in V,$ it suffices to verify the inequality \eqref{equivalent}  for all functions $f$ with $f(x)=0.$ Note that the inequality only involves the values of $f$ on $S_1(x).$ From now on, we label the vertices in $S_1(x)$ as $\{y_1,\cdots, y_M\}$ where $M=\dd_x.$ Any function $f$ on $S_1(x)$ can be understood as an $M$-tuple
$$(Y_1,\cdots, Y_M):=(f(y_1),\cdots, f(y_M)),$$ and the space of functions on $S_1(x)$ is identified with an $M$-dimensional vector space $\R^M$ indexed by the vertices of $S_1(x).$ 

For any $x,y\in V$ with $x\sim y,$ we denote \begin{equation}\label{quantalpha}\alpha_{xy}:=\frac{m_x}{\mu_{x{y}}}\left(\frac{\DD_x+\DD_{y}}{2}-\frac{2\mu_{xy}}{m_{y}}\right),\end{equation} which will be a key quantity in our argument, see the corollary below.

%To further simplify the notations, for any $y\in V,$ we write $Y:C(V)\to \R$ the evaluation function at $y$ as $$Y(f)=f(y),\quad \forall f\in C(V).$$ In the following, we use the capital Latin words for evaluation functions at corresponding vertices. Let 
%$\Omega=\{y_i\}_{1\leq i\leq M}\subset V,$ the function $f$ restricted to 
%$\Omega$ can be represented as
%$$(y_1,\cdots,y_M)\mapsto (Y_1(f),\cdots, Y_M(f)),$$ and $(Y_1,\cdots, Y_M)$ can be understood as $M$ variables labeled by $\{y_i\}_{1\leq i\leq M}.$

\begin{coro}\label{equCD0}If $\gir(x)\geq 5,$ then $\CCD(0,\infty, x)$ holds if and only if
\begin{equation*}\left(\sum_{y_i\sim x}\frac{\mu_{xy_i}}{m_x} Y_i\right)^2\geq \sum_{y_i\sim x} \frac{\mu_{x{y_i}}}{m_x}\left(\frac{\DD_x+\DD_{y_i}}{2}-\frac{2\mu_{xy_i}}{m_{y_i}}\right)Y_i^2,\quad\forall\ Y_i\in\R,\end{equation*} or equivalently,
\begin{equation}\label{char1}\left(\sum_{y_i\sim x} Y_i\right)^2\geq \sum_{y_i\sim x} \alpha_{xy_i}Y_i^2,\quad\forall\ Y_i\in\R,\end{equation} where $\alpha_{xy_i}$ is defined in \eqref{quantalpha}. 
\end{coro} 
\begin{proof} The first inequality is equivalent to \eqref{equivalent} for $K=0,n=\infty$. The second one follows from the first one by setting $Y_i'=\frac{\mu_{x{y_i}}}{m_x}Y_i$ for all $y_i\sim x$ and rename $Y_i'$ as $Y_i.$
\end{proof}

\begin{coro}\label{coropending} Let $x$ be a pending vertex,  i.e. $\dd_x=1,$ in a weighted graph $G.$ Then $\CCD(0,\infty,x)$
\begin{enumerate}\item always holds for normalized Laplacian, and
\item holds for unweighted physical Laplacian if and only if $\dd_y\leq 5,$ for $y\sim x.$
\end{enumerate}

\end{coro} \begin{proof}By the inequality \eqref{char1}, $\CCD(0,\infty,x)$ is equivalent to $$\frac{m_x}{m_y}\left(\frac12\frac{\mu_y}{\mu_x}-2\right)\leq \frac12,$$ where $y\sim x.$
\end{proof}

The following calculus lemma will be useful in our setting.
\begin{lemma}\label{calclem} Let $\{a_i\}_{1\leq i\leq N}, a_i\geq 0,$ and $c>0.$ The inequality, 
$$(Y+\sum_{i=1}^N Y_i)^2+\sum_{i=1}^N a_i Y_i^2\geq cY^2,$$ cannot hold for all $Y,Y_i\in \R$ $(1\leq i\leq N)$ if one of the following holds:
\begin{enumerate}
\item $c\geq 1.$
\item $c>0$ and $a_j=0$ for some $j\in\{1,\cdots,N\}.$
\end{enumerate}
\end{lemma}
\begin{proof} Suppose it holds for all $Y$ and $Y_i.$

\emph{(1)} $c\geq1.$ For any $t\in \R,$ setting $Y_i=t Y,$ $1\leq i\leq N,$ we have $$Y^2(1+Nt)^2+\sum_ia_it^2Y^2\geq cY^2.$$ This yields $2Nt+(N^2+\sum a_i)t^2\geq c-1\geq 0.$ It is not true for $$t\in\left(-2N(N^2+\sum a_i)^{-1},0\right).$$

\emph{(2)} $c>0$ and $a_j=0$ for some $j.$ Let $Y_i=0$ for all $i\neq j$ and $1\leq i\leq N.$ Then $$(Y+Y_j)^2\geq cY^2.$$ This yields a contradiction by setting $Y=-Y_j\neq 0.$
\end{proof}

\begin{lemma}\label{thm-strict} If $\gir(x)\geq 5$ and $\CCD(0,\infty, x)$ holds for a vertex $x$ in a weighted graph $G.$ Suppose that $\dd_x\geq 2,$ then $$ \alpha_{xy}<1,\quad \forall\ y\sim x.$$
\end{lemma}
\begin{proof} Without loss of generality, consider $y=y_1\sim x.$ Then setting $Y_i=0,$ for all $i\geq 2$ in \eqref{char1}, we have $$ \alpha_{xy_1}\leq1.$$ Suppose that $\alpha_{xy_1}=1.$ For the terms on the right hand side of \eqref{char1}, we eliminate those with positive coefficients on the right hand side, and move those with negative coefficients to the left hand side. This reduces to the case $(1)$ in Lemma~\ref{calclem}. This yields a contradiction and proves the strict inequality.
\end{proof}

For any $x\in V,$ we define $Q_x:=\{y\in V: y\sim x, \alpha_{xy}>0\}$ and $q_x:=\sharp Q_x.$ In fact the set $Q_x$ consists of those neighbors of $x$ which contribute positive terms on the right hand side of \eqref{char1}.

\begin{lemma}\label{lem1} If $\gir(x)\geq 5$ and $\CCD(0,\infty,x)$ holds for a vertex $x$ in a weighted graph $G.$ Then $q_x\leq 1.$
\end{lemma}
\begin{proof} Suppose that $q_x\geq 2.$ Without loss of generality, pick $y_1,y_2\in Q_x.$ By setting $Y_j=0$ for all $3\leq j\leq M$ in \eqref{char1}, we have
$$(Y_1+Y_2)^2\geq \alpha_{xy_1}Y_1^2+\alpha_{xy_2}Y_2^2\geq \alpha_{xy_1}Y_1^2, \quad \forall\ Y_1,Y_2\in\R.$$ This is impossible by $(2)$ in Lemma~\ref{calclem}.
\end{proof}

\section{Normalized Laplacians}\label{s:normal}
In this section, we consider the curvature dimension conditions for normalized Laplacians. For unweighted normalized Laplacians, a corollary of Theorem~\ref{main1} reads as follows.
\begin{coro} Let $G$ be an unweighted normalized graph, and for some $x\in V,$ $\gir(x)\geq 5.$ Then $\CCD(0,\infty,x)$ is equivalent to 
\begin{equation}\label{char2}\left(\sum_{y_i\sim x} Y_i\right)^2\geq \dd_x\sum_{y_i\sim x}\left(1-\frac{2}{\dd_{y_i}}\right)Y_i^2,\quad\forall\ Y_i\in\R.\end{equation}
\end{coro}

In this setting, $\alpha_{xy}=\dd_x\left(1-\frac{2}{\dd_{y}}\right)$ for any $x\sim y$ and $Q_x:=\{y\sim x: \dd_y\geq 3\}$ for all $x\in V.$ By Lemma~\ref{lem1}, we know that $q_x\leq 1$ if $\gir(x)\geq 5$ and $\CCD(0,\infty,x)$ hold.
%For any $x\in V,$ we define $Q_x:=\{y\sim x: \dd_y\geq 3\}$ and $q_x:=\sharp Q_x.$ In fact the set $Q_x$ consists of those neighbors of $x$ which contribute positive terms on the right hand side of \eqref{char2}.
%We list neighbors of $x$ as $\{y_1,\cdots,y_M\}$ where $M=\dd_x$ in the following.

%\begin{lemma}\label{lem1} Let $(V,E,m,\mu)$ be equipped with an unweighted normalized Laplacian, and $\gir(x)\geq 5$ and $\CCD(0,\infty,x)$ hold for some $x\in V.$ Then $q_x\leq 1.$
%\end{lemma}
%\begin{proof} Suppose that $q_x\geq 2.$ Without loss of generality, we set $y_1,y_2\in Q_x.$ By setting $Y_j=0$ for all $3\leq j\leq M$ in \eqref{char2}, we have
%$$(Y_1+Y_2)^2\geq c_1Y_1^2+c_2Y_2^2,\quad \forall\ Y_1,Y_2\in\R$$ where $c_i=\dd_x(1-2/\dd_{y_i})>0,\ i=1,2.$ This is impossible by choosing $Y_1=-Y_2\neq 0.$ A contradiction.
%\end{proof}

\begin{lemma}\label{lem2} Let $G$ be an unweighted normalized graph and for some $x\in V$ with $\dd_x\geq 2,$ $\gir(x)\geq 5$ and $\CCD(0,\infty,x)$ hold. If $q_x=1,$ then $\dd_x=2,$ $\dd_{y_1}=3$ and $\dd_{y_2}=1$ where $y_i\sim x, i=1,2.$
\end{lemma}
\begin{proof} Without loss of generality, let $y_1\in Q_x,$ i.e. $\dd_{y_1}\geq 3,$ and $S_1(x)=\{y_1,\cdots,y_M\}$ with $M=\dd_x,$ where $\dd_{y_{i}}\leq 2,$ for all $i\neq 1.$ By Lemma~\ref{thm-strict}, $\alpha_{xy_1}<1.$ By this inequality,
$$\dd_{y_1}<2\left(1-\frac{1}{\dd_x}\right)^{-1}\leq 4,$$ which yields that $\dd_{y_1}=3,$ and
$$\dd_x<\left(1-\frac{2}{\dd_{y_1}}\right)^{-1}=3,$$ which implies $\dd_x=2.$ Hence $S_1(x)=\{y_1,y_2\}.$

For the case of  $\dd_{y_2}=2,$ by \eqref{char2},
$$(Y_1+Y_2)^2\geq \frac{2}{3}Y_1^2,\quad \forall\ Y_1,Y_2\in\R,$$ which is impossible by setting $Y_1=-Y_2\neq 0,$ or $(2)$ in Lemma~\ref{calclem}.

For the case of $\dd_{y_2}=1,$ by \eqref{char2},
$$(Y_1+Y_2)^2+2Y_2^2\geq \frac{2}{3}Y_1^2,\quad \forall\ Y_1,Y_2\in\R,$$
which is true and proves the lemma.
\end{proof}

Combining Corollary~\ref{coropending}, Lemma~\ref{lem1} with Lemma~\ref{lem2}, we have the following.
\begin{lemma}\label{lem3} Let $G$ be equipped with an unweighted normalized Laplacian. For $x\in V,$ if $\gir(x)\geq 5$ and $\CCD(0,\infty,x)$ holds, then we have three cases:
\begin{enumerate}\item $\dd_x=1,$ then $\dd_y$ is arbitrary for $y\sim x.$

\item $\dd_x=2,$ then $(a)$ either $\dd_{y_1}=3$ and $\dd_{y_2}=1,$ or $(b)$ $\dd_{y_1}\leq 2$ and $\dd_{y_2}\leq 2,$ where $y_i\sim x, i=1,2.$ 

\item $\dd_x\geq 3,$ then $\dd_{y}\leq 2$ for all $y\sim x.$ \end{enumerate}
\end{lemma}
\begin{proof} The case $(1)$ follows from Corollary~\ref{coropending}. For $q_x=0,$ i.e. $\dd_{y}\leq 2$ for all $y\sim x,$ \eqref{char2} always holds. For $q_x=1,$ we apply Lemma~\ref{lem2}. Hence, we obtain the cases $(2)$ and $(3).$ 
\end{proof}

Now we are ready to classify the graph with large girth which has nonnegative curvature dimension condition. The first case is that the degree of all vertices are at least two.
\begin{proof}[Proof of Theorem~\ref{main-normalized1}] Applying Lemma~\ref{lem3} at the vertex $x_0,$ noting that $\inf_{x\in V}\dd_x\geq 2,$ we have $\dd_{x_0}\geq 2,\dd_{y}=2$ for all $y\sim x.$ 

We claim that $\dd_{x_0}=2.$ 
Suppose not, i.e. $\dd_{x_0}\geq 3.$ Pick a neighbor of ${x_0},$ say $y_1.$ Noting that $\dd_{y_1}=2$ and $\gir(x_0)\geq5,$ we have $\gir(y_1)\geq 5.$ Now we may apply Lemma~\ref{lem3} to $y_1,$ and obtain that $\dd_{x_0}=2.$ A contradiction. 

Hence $\dd_{x_0}=2,\dd_{y}=2$ for all $y\sim x.$ This yields that $\gir(y)\geq 5$ for $y\sim x.$ Using the same argument at $y,$ we have that $\dd_z=2$ for $z\sim y.$ Continuing this process, we conclude that $\dd_x=2,\forall x\in V.$ This proves the theorem.

\end{proof}

%Note that in this theorem, we only need to assume that the graph has a vertex whose girth is at least $5$ if all vertices has degree at least $2.$ 

Next we classify the general cases without any restrictions on vertex degrees.

\begin{proof}[Proof of Theorem~\ref{main-normalized2}] We claim that $W_3:=\sharp\{x\in V:\dd_x\geq 3\}\leq 1.$
Suppose not. Let $x_1,x_2$ be two distinct vertices with $\dd_{x_i}\geq 3, i=1,2.$ Then by the connectedness, there is a walk connecting them, $x_1=z_0\sim z_1\sim \cdots\sim z_N=x_2$ with $N\geq 1.$ Applying Lemma~\ref{lem3} at the vertex $x_1,$ we have $\dd_{z_1}\leq 2$ and $N\geq 2,$ which implies that $\dd_{z_1}=2$ since $z_1$ lies in the walk connecting $x_1$ and $x_2.$ Applying Lemma~\ref{lem3} at the vertex $z_1,$ we obtain that $\dd_{x_1}=3$ and $\dd_{z_2}=1.$ This yields a contradiction since $z_2$ lies in the walk and proves the claim.

For the case of $W_3=0,$ we get the classification $(a)$ in the theorem.

Let $W_3=1$ and $\dd_{x_0}\geq 3$ for $x_0\in V.$ By Lemma~\ref{lem3}, $\dd_y\leq 2$ for any $y\sim x_0.$ We divide it into cases:

\emph{Case 1.} $\dd_{x_0}\geq 4.$ Applying Lemma~\ref{lem3} at any $y\sim x,$ we have $\dd_y=1.$ Hence the graph is $\str_n,$ $n\geq 4.$

\emph{Case 2.} $\dd_{x_0}=3.$ For any $y\sim x_0$ satisfying $\dd_y=2,$ applying Lemma~\ref{lem3} to $y,$ we get $\dd_z=1$ for $z\sim y$ $(z\neq x).$ 
Hence we obtain $\str_3,$ or $\str_3^{i},$ $i=1,2,3.$ This gives the classification $(b)$.
\end{proof}

\section{Physical Laplacians}\label{s:physical}
In this section, we consider unweighted physical Laplacians and have a corollary of Theorem~\ref{main1} as follows.
\begin{coro} Let $G$ be an unweighted physical graph and for $x\in V$ $\gir(x)\geq 5.$  Then $\CCD(0,\infty,x)$ is equivalent to 
\begin{equation}\label{char3}\left(\sum_{y_i\sim x} Y_i\right)^2\geq \sum_{y_i\sim x}\left(\frac{\dd_x+\dd_y}{2}-2\right)Y_i^2,\quad\forall\ Y_i\in\R.\end{equation}
\end{coro}

In this setting, $\alpha_{xy}=\frac{\dd_x+\dd_y}{2}-2$ for any $x\sim y$ and $Q_x:=\{y\in V: y\sim x, \dd_x+\dd_y\geq 5\}$ for all $x\in V.$ 
%Moreover, by Lemma~\ref{lem1}, $q_x\leq 1$ if $\gir(x)\geq 5$ and $\CCD(0,\infty)$ hold.

\begin{lemma}\label{lem6} Let $G$ be an unweighted physical graph and for $x\in V,$ $\gir(x)\geq 5.$ Then $\CCD(0,\infty,x)$ holds if and only if we have the following cases:
\begin{enumerate}\item $\dd_x=1,$ then $\dd_y\leq 5$ for $y\sim x.$

\item $\dd_x=2,$ then $\dd_{y_i}\leq 2$ for $y_i\sim x, i=1,2.$ 

\item $\dd_x\geq 3,$ then $\dd_{y_i}=1$ for all $y_i\sim x, 1\leq i\leq 3.$ \end{enumerate}
\end{lemma}
\begin{proof} Without loss of generality, we may assume $\dd_x\geq 2$ by Corollary~\ref{coropending} and denote $S_1(x)=\{y_1,\cdots,y_M\}$ with $M=\dd_x.$ By Lemma~\ref{lem1}, $q_x\leq 1.$ We divide it into cases:

\emph{Case 1.} $q_x=1.$ Without loss of generality, let $y_1\in Q_x.$ By Lemma~\ref{thm-strict}, we have
$$5\leq \dd_x+\dd_{y_1}<6,$$ which yields $\dd_x+\dd_{y_1}=5.$ Then $\dd_x=2,3$ or $4.$ The case $\dd_x=4$ can be excluded by $\dd_x+\dd_{y_2}\leq 4$ since $q_x=1.$ 

Suppose $\dd_x=2.$ Then $\dd_{y_1}=3,\dd_{y_2}\leq2.$ For $\dd_{y_2}=2,$ the equation \eqref{char3} reads as $$(Y_1+Y_2)^2\geq \frac12 Y_1^2,\quad \forall\ Y_1,Y_2\in \R,$$ which is impossible.
For $\dd_{y_2}=1,$ we have $$(Y_1+Y_2)^2\geq \frac12Y_1^2-\frac12Y_2^2,\quad \forall\ Y_1,Y_2\in \R,,$$ which is also wrong.

Suppose $\dd_x=3.$ Then $\dd_{y_1}=2,\dd_{y_2}=\dd_{y_3}=1.$ Then the equation \eqref{char3} has the form
$$(Y_1+Y_2+Y_3)^2\geq \frac12 Y_1^2,\quad\forall\ Y_1,Y_2,Y_3\in\R,$$ which is excluded by $(2)$ in Lemma~\ref{calclem}.

\emph{Case 2.} $q_x=0.$ That is, $\dd_x+\dd_{y_i}\leq 4$ for all $y_i\sim x.$ This implies that $\dd_x\leq 3.$
For the case of $\dd_x=2,$ we have $\dd_{y_i}\leq 2,$ for $i=1,2.$ This gives the case (2) in the lemma.
For the case of $\dd_x=3,$ $\dd_{y_i}=1,$ for $1\leq i\leq 3.$ That is the case $(3)$.

This proves the lemma.
\end{proof}

By this lemma, following the arguments as in Theorem~\ref{main-normalized1} and \ref{main-normalized2}, one can prove the following results for unweighted physical Laplacians. We omit the proofs here.

\begin{thm}Let $G$ be an unweighted physical graph satisfying the $\CCD(0,\infty)$ condition and $\inf_{x\in V}\dd_x\geq 2.$ If for some $x_0\in V,$ $\gir(x_0)\geq 5,$ then the graph is either the infinite line $P_{\Z}$ or the cycle graphs $C_n$ for $n\geq 5.$
\end{thm}

\begin{thm} Let $G$ be an unweighted physical graph satisfying the $\CCD(0,\infty)$ condition with girth at least five. Then the graph is one of the following:
\begin{enumerate}
\item $P_k$ $(k\geq 1),$ $C_n$ $(n\geq 5),$ $P_{\Z},$ or $P_{\N}.$
\item $\str_n$ $(n\geq 3).$ 
\end{enumerate}
\end{thm}

{\bf Acknowledgements.} 
We thank Dr. Yan Huang for drawing the pictures in the paper.

B. H. is supported by NSFC, grant no. 11401106. Y. L. is supported by NSFC, grant no. 11271011, the Fundamental Research Funds for the Central Universities and the Research Funds of Renmin University of China($11$XNI$004$).
\bibliography{girth}
\bibliographystyle{alpha}

\end{document}